\documentclass[12pt]{amsart}

\usepackage{graphicx}
\usepackage{amsthm}
\usepackage{amsmath}
\usepackage{amssymb}
\usepackage{color}
\usepackage{diagrams}
\usepackage{hyperref}
\usepackage{multirow}
\usepackage{makecell}
\usepackage{mathtools}
\usepackage{comment}
\usepackage{tikz-cd}	
\usepackage{mathrsfs}   
\usepackage{float}

\usepackage{setspace}

\newtheorem{Thm}{Theorem}[section]
\newtheorem{Lem}[Thm]{Lemma}
\newtheorem{Cor}{Corollary}[Thm]
\newtheorem{Prop}[Thm]{Proposition}

\theoremstyle{definition}
\newtheorem{Def}[Thm]{Definition}
\newtheorem{eg}[Thm]{Example}

\theoremstyle{remark}
\newtheorem{Rmk}[Thm]{Remark}

\newenvironment{pf}{\begin{proof}}{\end{proof}}

\numberwithin{equation}{section}

\newcommand{\Z}{\mathbf{Z}}
\newcommand{\Q}{\mathbf{Q}}

\newcommand{\C}{\mathbf{C}}
\newcommand{\F}{\mathbf{F}}

\newcommand{\RNum}[1]{\uppercase\expandafter{\romannumeral #1\relax}}
\def\lrd{\overwithdelims()}

\DeclareFontFamily{U}{wncy}{}
\DeclareFontShape{U}{wncy}{m}{n}{<->wncyr10}{}
\DeclareSymbolFont{mcy}{U}{wncy}{m}{n}
\DeclareMathSymbol{\Sha}{\mathord}{mcy}{"58} 

\begin{document}
\title{Galois action on the N\'eron-Severi group of Dwork surfaces}

\author{Lian Duan}
\address{Department of Mathematics and Statistics,  University of Massachusetts Amherst, 710 N. Pleasant Street Amherst, MA 01003-9305, USA}
\curraddr{Department of Mathematics and Statistics,  University of Massachusetts Amherst, 710 N. Pleasant Street Amherst, MA 01003-9305, USA}
\email{duan@math.umass.edu}






\begin{abstract}
We study the Galois action attached to the Dwrok surfaces $X_{\lambda}:X_0^4+X_1^4+X_2^4+X_3^4-4\lambda X_0X_1X_2X_3=0$ with parameter $\lambda$ in a number field $F$. We show that when $X_{\lambda}$ has geometric Picard number $19$, its N\'eron-Severi group $NS(\overline{X}_{\lambda})\otimes \Q$ is a direct sum of quadratic characters. We provide two proofs to this conclusion in our article. In particular, the geometrically proof determines the conductor of each of quadratic characters. Our result matches the one in \cite{Voight2}. With this decomposition, we give another proof to a result of Wan \cite{Wan-Mirror-Sym}.
	
\end{abstract}

\maketitle
\tableofcontents
	\section{Introduction}
	Let $F$ be a number field. For a fixed integer $n$, the family of parametric hypersurfaces defined by 
	$$
	X_0^n+X_2^n+\cdots +X_{n-1}^n=\lambda X_0X_1\cdots X_{n-1}
	$$
	is called the Dwork family of degree $n$. This family appears in Dwork's articles (\cite[page. 249]{Dwork-Deformation}, \cite[$\S $6.25, 6.30]{Dwork-p-cycle}) as examples of his deformation theory. Since then, it has been studied in many areas of mathematics, cf.~\cite{Katz-look-at-Dwork}. In particular, the zeta functions of this family have been studied, for example, by Wan \cite{Wan-Mirror-Sym} and Goutet \cite{Goutet-Dwork}.
	In this paper, we focus on the case $n=4$, and modify the defining equation of the family as follows:
	$$
	X_{\lambda}: X_0^4+X_1^4+X_2^4+X_3^4-4\lambda X_0X_1X_2X_3=0
	$$
	For each $\lambda\in F$ such that $\lambda^4\neq 1$, we get a smooth and hence $K3$ projective surface $X_{\lambda}$. Every smooth $X_{\lambda}$ has trivial first and third singular cohomology and torsion free second cohomology with rank $22$. If we denote by $\overline{X}_{\lambda}=X_{\lambda}\times_{F}\overline{F}$ the base change of $X_{\lambda}$ to the algebraic closure of $F$, then the geometric Picard group $Pic(\overline{X}_{\lambda})$ is equal to the geometric Neron Severi group $NS(\overline{X}_{\lambda})$, which is a free $\Z$-module of rank $19$ (nonsigular) or $20$ (singular) and in general $19$, (\cite[remark 4.5]{Bini-Garbagnati}, \cite{Elkies-Schutts-K3}). As a set of equivalence classes of algebraic cycles, $NS(\overline{X}_{\lambda})$ admits an action by $G_F:=Gal(\overline{F}/F)$. Let ${a \lrd \bullet}$ be the quadratic character corresponding to $F(\sqrt{a})$, where $a$ a square free algebraic integer in $F$. In this paper we prove the following theorem.
	
	\begin{Thm}\label{new main Thm}
	If $X_{\lambda}$ is smooth with Picard rank $19$, then as a $G_{F}$-representation space, $NS(\overline{X}_{\lambda})\otimes \Q$ is a direct sum of 1-dimensional linear characters. Specifically, 
		\begin{equation*}
			\resizebox{1.0 \textwidth}{!} 
			{
			$NS(\overline{X}_{\lambda})\otimes \Q=id\bigoplus {-(\lambda^2-1) \lrd \bullet}^{\oplus 3} \bigoplus {-(\lambda^2+1) \lrd \bullet}^{\oplus 3}\bigoplus {-2(\lambda^4-1) \lrd \bullet}^{\oplus 6}\bigoplus {2(\lambda^4-1) \lrd \bullet}^{\oplus 6}.$
			}
			\end{equation*}
	\end{Thm}

	When $\lambda=0$, we can write down the algebraic cycles explicitly (\cite[$\S$ 4]{Bini-Garbagnati}, \cite[chapter 17 $\S 1$]{K3-lecture-notes}), which then allows us to compute the Galois action on $NS(\overline{X}_{\lambda})$ via these algebraic cycles. However, for general $\lambda$, we have not found a concrete description of the corresponding result. In our paper, we study the Galois action on $NS(\overline{X}_{\lambda})$ in two ways. 
	
	In the first method, note that for a general $\lambda$ the surface $X_{\lambda}$ has a large group of $\Q$-rational automorphisms. The action of these automorphisms on $X_{\lambda}$ induces an action on $NS(\overline{X}_{\lambda})\otimes_{\Z}\Q$ that commutes with the $G_F$-action. Thus we can use these automorphisms to decompose $NS(\overline{X}_{\lambda})\otimes_{\Z} \Q$ into a direct sum of $G_F$-submodules. Furthermore,  by replacing the whole automorphism group with its subgroups, we are able to decompose every direct summand into smaller $G_F$-modules. By repeating this process we are end up with nineteen $1$-dimensional $G_F$-representations. Then by Poincare duality, each $1$-dimensional $G_F$-representation corresponds to a quadratic character. Then, by comparing the zeta functions of $X_{\lambda}$ with that of one of its Mirror symmetry, we prove theorem \ref{new main Thm} but without specifying the conductors. For more details, see Theorem \ref{main Thm}. We then computed $NS(\overline{X}_{\lambda})\times \Q$ for specific values of $\lambda\in \Q$ by counting points. These experimental data led us to conjecture the precise form of theorem \ref{new main Thm} and verified it in special cases.
	
	In a series of work, Doran, Kelly, Salerno, Sperber, Voight and Whitcher (\cite{Voight1}, \cite{Voight2}) carry out a detailed study of the Picard-Fuchs equations and the zeta functions of several pencils of K3 surfaces or geometric Picard number $19$. They explicitly determine the factors of the L-functions of these pencils in terms of hypergeometric sums. In particular, one of their pencils ($F_4$ in their notation) is our Dwork family, upon learning our conjecture, they apply their method to give a complete proof of Theorem \ref{new main Thm} (cf. \cite[Theorem 1.4.1 (a) and expressions (1.4.6)]{Voight2}).	
	
	Remke Kloosterman and John Voight also kindly suggest that we apply the method of \cite{Kloosterman} to give a geometric proof to Theorem \ref{new main Thm}. Specifically, we consider several degree $2$ quotients of $X_{\lambda}$ which are induced by involutions of $X_{\lambda}$. Each of these quotients is not only a degree $2$ del Pezzo surface, but also a double covering of $\mathbb{P}^2$ branched along a quartic curve. By pulling back the $56$ lines lying above the bitangents of the quartic curve, we can write down explicit generators of $NS(X_{\lambda})\otimes \Q$, and hence completely prove Theorem \ref{new main Thm}.
	
	Note that the explicit generators found by the second method provide more information than the decomposition in theorem \ref{new main Thm} and will be useful for other problems. We are in the process of carrying out this analysis for the other families in \cite{Voight2}. More generally, our initial representation theory method, while not sufficient to completely prove our main theorem, is more flexible and can be applied to other situations. We will pursue these and related ideas in our future work. 
		
	As an application of our main theorem, we consider the family of surfaces $Y_{\lambda}$ coming from the resolution of the singularities of the hypersurface $M_{\lambda}$ below: 
	$$
	M_{\lambda}:(Y_0+Y_1+Y_2+Y_3)^4-(4\lambda)^4Y_0Y_1Y_2Y_3=0
	$$ 
	If $\lambda^4\neq 1$ and $\lambda\neq 0$, then both $X_{\lambda}$ and $Y_{\lambda}$ are smooth $K3$ surfaces and have the same Picard number (\cite{Elkies-Schutts-K3} or \cite[$\S 3.2.2$]{Bini-Garbagnati}). Wan (\cite{Wan-Mirror-Sym}) uses a Gauss sum calculation to deduce a congruence relation of the numbers of rational points over finite fields between the Dwork family of dimension $d$ and their strong mirror pair. When the dimension is $2$, such pair is our $X_{\lambda}$ and $Y_{\lambda}$. By studying the transcendental parts of $X_{\lambda}$ and $Y_{\lambda}$, we give another proof of Wan's result for surface case. Note that to deduce this corollary, argument only makes use of the decomposition of $NS(\overline{X}_{\lambda})\otimes \Q$ in the form of theorem \ref{main Thm} and not the explicit determination of the conductor of the characters.
	\begin{Cor}\cite[Thm.~1.1]{Wan-Mirror-Sym}\label{Cor}
		With $X_{\lambda}$ and $Y_{\lambda}$ defined above, we have 
		$$
		\# X_{\lambda}(\F_{p^k})\equiv \# Y_{\lambda}(\F_{p^k})\ (mod\ p^k)
		$$
	\end{Cor}

	We now give an outline of our paper. In section 2 we review some basic facts about $K3$ surfaces. In Section 3 we consider the group $H:=S_4\ltimes (\Z/2\Z)^{2}$, note that $H$ has natural group action on $X_{\lambda}$, which is commutative with the Galois action. Thus we can decompose $H^2_{et}(\overline{X}_{\lambda}, \overline{\Q}_{\ell}(1))$ canonically as $H$-representation. After excluding the transcendental part $T_{\lambda}$ and using Poincare duality, we deduce most parts of our Theorem \ref{main Thm} (except parts (b), (d)). 
	
	In Section 4, we consider another family of $K3$ surfaces $Y_{\lambda}$. By comparing our situation with the classical double covering case, we prove that if $\lambda^4\neq 1$, the transcendental parts of the pair $X_{\lambda}, Y_{\lambda}$ have the same $G_F$-representation. As a consequence we finish the proof of our main theorem. In addition, we prove corollary \ref{Cor}. 
	
	In Section 5, we review the necessary properties of del Pezzo surfaces, especially for those of degree $2$. Then we focus on the double quotients of $X_{\lambda}$, and use deformation theory to find out the generators of their N\'eron-Severi groups. 
	
	In Section 6, we compute the Galois action on $NS(\overline{X}_{\lambda})$ based on the generators found in section 5, and thus give a complete proof to our main theorem.

	To simplify our arguments, we will adopt the following notations. In this paper, unless otherwise mentioned specifically, all varieties are algebraic, smooth, projective and defined over number field. Given a variety $X$ defined over $F$, we denote by $\overline{X}:=X\times_{F}\overline{F}$ the base change of $X$ to the algebraic closure of the ground field. The notation $X_{\lambda}$ is reserved to refer to the elements in Dwork family with parameter $\lambda$. 
	
	Given a number field $K$ and its Galois closure $\overline{K}$, $G_K$ refers to the absolute Galois group $Gal(\overline{K}/K)$. In this paper, we will usually switch between complex manifold cohomology and \'{e}tale cohomology. To differentiate them, for a smooth variety $X$, we use $H^i$ without subindex to refer the former, and $H^i_{et}$ for the latter. We will abuse the notation $NS(\overline{X}_{\lambda})$ to refer both to its image in $H^2_{}(\overline{X}_{\lambda},\Z)$ and $H^2_{et}(\overline{X}_{\lambda},\Z_{\ell}(1))$ when there is no confusion. Also, we denote by $T_{\lambda}(\overline{X})$ the transcendental lattice of $\overline{X}$ in both $H^2_{}(\overline{X},\Z)$ and $H^2_{et}(\overline{X}, \Z_{\ell}(1))$.  	
	
	\subsection*{Acknowledgement}
	We would like to thank Professor Siman Wong for suggesting this interesting topic, and giving helpful advice. We thank Professor Remke Kloosterman and Professor John Voight for sharing their beautiful ideas with the author and to Professor Paul Hacking and Professor Eyal Markman for conversation and suggestions in algebraic geometry.

	\section{Background}
	
	This section is intended as a review of the background needed for the rest of the paper. We list the necessary definitions and properties. We usually skip the proofs, but cite the references for interested readers. All of the definitions and basic properties used can be found in many textbooks on algebraic surfaces and \'{e}tale cohomology. In this paper, for the geometry of $K3$ surfaces, we follow the first two chapters of \cite{K3-lecture-notes}, for the \'{e}tale cohomology background, we follow \cite{MilEC}. 
	\begin{Def}
		A $K3$ surface over field $k$ is a complete, non-singular surface $X$ such that 
		$$
		\Omega^2_{X/k}\simeq \mathcal{O}_X \text{ and } H^1(X,\mathcal{O}_X)=0
		$$
	\end{Def}
	
	\begin{Prop}\cite[example 1.1.3]{K3-lecture-notes}
		Any smooth quartic surface in $\mathbb{P}^3$ is $K3$.
	\end{Prop}
	
	For any smooth surface $X$, we denote by $Pic(\overline{X})$ the geometric Picard group which is generated by linear equivalent classes of algebraic cycles of $\overline{X}$,  and we denote by $Num(\overline{X})$ the group generated by numerical equivalent classes of divisors.

	\begin{Def}
		The geometric Neron-Severi group of an algebraic surface $X$ is the quotient 
		$$
		NS(\overline{X}):=Pic(\overline{X})/Pic^0(\overline{X})
		$$
		where $Pic^0(\overline{X})$ is the subgroup of $Pic(\overline{X})$ of line bundles that are algebraically equivalent to the zero divisor. In other words, $NS(\overline{X})$ is the group of algebraically equivalent classes. 
	\end{Def}
	
	\begin{Prop}\cite[prop.~1.2.4]{K3-lecture-notes}
		For a $K3$ surface $X$, 
		$$
		Pic(\overline{X})\simeq NS(\overline{X}) \simeq Num(\overline{X})
		$$
	\end{Prop}
	
	We call the rank of $NS(\overline{X})$ the (geometric) Picard number of $X$ and denoted it by $\rho(\overline{X})$. 
	
	\begin{Prop}[Hodge index theorem] \cite[chap.~1, $\S$ 2.2]{K3-lecture-notes}\label{ind thm}
		The signature of the intersection form on $Num(\overline{X})$ for a smooth surface is $(1, \rho(\overline{X})-1)$. 
	\end{Prop}
	
	\begin{Cor}
		When $X$ is $K3$, the intersection form on $NS(\overline{X})$ is nondegenerate and has signature $(1, \rho(\overline{X})-1)$.
	\end{Cor}

	Many results about $K3$ surfaces are related to their cohomological properties. In this paper, we will consider two kinds of cohomology, the singular cohomology $H^i(\overline{X}, \Z)$ and the \'{e}tale cohomology $H^i_{et}(\overline{X}, \Z_{\ell}(r))$, where $r$ stands for the $r$-th Tate twist. By the comparison theorem \cite[chap.~\RNum{3}, thm.~3.12]{MilEC}, we know that for $K3$ surfaces, both of the cohomology are of the same rank. In particular, $H^2(\overline{X}, \Z)$ is a free $\Z$-module of rank $22$, hence so is $H^2_{et}(\overline{X}, \Z_{\ell})$ as a $\Z_{\ell}$-module. 
	As $NS(\overline{X})$ embeds into both of the two cohomology groups, we have:
	\begin{Cor}\label{signature}
		The geometric Neron-Severi group of a $K3$ surface is a free module.
	\end{Cor}
	By the Lefschetz theorem on (1,1) classes we know that $NS(\overline{X})$ embeds into $H^{1,1}(\overline{X},\C)$. Since $h^{1,1}(\overline{X})=20$ \cite[chap.~1, $\S$ 2]{K3-lecture-notes}, it follows that $\rho({\overline{X}})\leq 20$. If a $K3$ surface is defined over a field of characteristic $0$ and has Picard number $20$, we call it a singular $K3$ surface (although it is still smooth). Otherwise, we say it is nonsingular. 
	
	\begin{Rmk}
		In this paper we will focus on $K3$ surfaces defined over number fields. If the characteristic of the ground field for the $K3$ surface is positive, then $\rho$ could be 22. 
	\end{Rmk}
	
	\begin{Prop}\cite[remark.~4.5]{Bini-Garbagnati} or \cite{Elkies-Schutts-K3}.
		For a Dwork surface $X_{\lambda}$, we have $19\leq \rho(\overline{X}_{\lambda})\leq 20$, and for general $\lambda$ we have $\rho(\overline{X}_{\lambda})=19$.
	\end{Prop}

	Poincar\'e duality gives rise to well-defined bilinear forms on $H^2_{}(\overline{X},\Z)$ and $H^2_{et}(\overline{X}, \Z_{\ell}(1))$. In this paper, we will use the \'{e}tale cohomology version. 
	
	\begin{Prop}\cite[chap.~\RNum{6}, cor.~11.2]{MilEC}\label{PD}
		Given a surface $X$,  the cup-product pairing
		$$
		H^2_{et}(X, \Z_{\ell}(1))\times H^2_{et}(X, \Z_{\ell}(1))\to H^4_{et}(X, \Z_{\ell}(2))
		$$
		is nondegenerate.
	\end{Prop}
	With this bilinear form we have a lattice structure (that is, a free $\Z$ or $\Z_{\ell}$ module with a bilinear form \cite[chap.~14]{K3-lecture-notes}, or \cite{Morrison-K3-SI}). A sublattice of a given lattice is a free submodule with inherited bilinear form. If $L'$ is a sub lattice of $L$, we call $L'$ primitive if $L/L'$ is torsion free. By the Lefschetz theorem on (1,1) class, we can view $NS(\overline{X})$ as a sublattice of $H^2(\overline{X},\Z)$ or $H^2_{et}(\overline{X}, \Z_{\ell}(1))$.
	
	\begin{Def}
		Under the above structure, the primitive sublattice of $H^2_{}(\overline{X},\Z)$, which contains $H^{2,0}(X,\C)$ and $H^{0,2}(X,\C)$ after base extension to complex field $\C$, is called the transcendental lattice, denoted by $T(X)$ or $T$ for simplicity. 
	\end{Def}
	
	\begin{Prop}\cite[Lem.~3.3.1]{K3-lecture-notes}.
		The transcendental lattice of a complex $K3$ surface is the orthogonal complement of the Neron-Severi group in $H^2(\overline{X},\Z)$, i.e.
		$$
		T(X)=NS(\overline{X})^{\bot}
		$$
	\end{Prop}

	\begin{Rmk}\label{Gal compatible}
		Suppose $X$ is a $K3$ surface defined over a number field $F$. Then both $NS(\overline{X})$ and $ T(X)$ are $G_F$-modules. Poincar\'e duality is compatible with $G_F$-action. Moreover, 
		$$
		NS(\overline{X})\otimes_{\Z} \Z_{\ell}\hookrightarrow H^2_{et}(\overline{X}, \Z_{\ell}(1))
		$$
		is an embedding of $G_F$-modules \cite[chap.~\RNum{4}, remark 9.6]{MilEC}. 
	\end{Rmk}
	
	To simplify our notations, we will view $NS(\overline{X})$ and $T(X)$ as sublattices of $H^2_{}(\overline{X},\Z)$ or $H^2_{et}(\overline{X},{\Z}_{\ell}(1))$ or $H^2_{et}(\overline{X},\overline{\Q}_{\ell}(1))$ when there is no confusion.
	
	We end this section by reviewing the canonical decomposition of the representation of a finite group \cite[$\S$.~2.6]{Serre-LRF}. If $H$ is a finite group and $\rho: H\to Aut(V)$ is a finite dimensional vector space over an algebraically closed field $k$ of characteristic $0$, then character theory \cite[$\S$.~2.6]{Serre-LRF} gives rise to a decomposition of $V$ as an $H$-representation which is independent of the choice of basis:
	$$
	V=\bigoplus\limits_{\chi_i} im\left( \frac{n_i}{|H|}\sum\limits_{h\in H}\chi_i(h^{-1})\rho(h)\right) ,
	$$ 
	where each $\chi_i$ is an irreducible character of $H$ and $n_i$ is the degree of $\chi_i$. We call the above decomposition the canonical decomposition of $V$ as an $H$-representation. Notice that each direct summand only depends on its corresponding character, thus we will denote by $V_{\chi_i}$ the direct summand corresponding to $\chi_i$. Now suppose that $V$ also admits the action from a compact group $G$. If the actions of $G$ and $H$ commute, i.e. $gh(x)=hg(x)$ for all $g\in G$, $h\in H$ and $x\in V$, then we have the following proposition.
	
	\begin{Prop}\label{Canonical dec}
		Given groups $G$ and $H$ and vector space $V$ as above, then every $V_{\chi_i}$ is also a sub $G$-representation.
	\end{Prop}
	
	\begin{pf}
		It is sufficient to show that every $V_{\chi_i}$ is stable under the action of $G$. To see this, assume $y\in V_{\chi_i}$. Then we can find $x\in V$ such that $y=\frac{n_i}{|H|}\sum\limits_{h\in H}\chi_i(h^{-1})h(x)$. Now for $g\in G$, we have 
		\begin{align*}
		g(y)=g\left( \frac{n_i}{|H|}\sum\limits_{h\in H}{\chi_i}(h^{-1})h(x)\right)&=\frac{n_i}{|H|}\sum\limits_{h\in H}{\chi_i}(h^{-1})g(h(x))\\
		&=\frac{n_i}{|H|}\sum\limits_{h\in H}{\chi_i}(h^{-1})h(g(x))\in V_{\chi_i}.
		\end{align*}

		Hence, $V_{\chi_i}$ is stable under the action of $G$.
	\end{pf}

	\section{The Geometric Group Action}
	Our goal in this and the next sections are to prove the weaker version of the main theorem.
	
	\begin{Thm}\label{main Thm}
					If $X_{\lambda}$ is smooth with Picard rank $19$, then as $G_F$-representation spaces, $NS(\overline{X}_{\lambda})\otimes \Q$ is a direct sum of 1-dimensional linear characters. Specifically, 
					\begin{equation}\label{Eqn-main}
					NS(X_{\lambda})\otimes \Q=id \oplus{a_1 \lrd \bullet}^{\oplus n_1}\oplus\cdots \oplus {a_r \lrd \bullet}^{\oplus n_r},
					\end{equation}
					such that:
					\begin{enumerate}
						\item[(a)] The trivial character $id$ comes from the hyperplane section;
						\item[(b)] Every $n_i$ is divisible by $3$;
						\item[(c)] At least one of the $n_i\geq 6$;
						\item[(d)] The index number $r\leq 5$.
					\end{enumerate}
				\end{Thm}
			
	In this section, we prove the decomposition in expression (\ref{Eqn-main}), (a), (c) and part of (b) of Theorem \ref{main Thm}. First we state out idea. Note that by the symmetry of the defining equation, each $X_{\lambda}$ admits an action from the group $H:=S_4\ltimes (\Z/2\Z)^2$ with the following generators:
	if $\sigma\in S_4$, then
	$$
	\sigma:(X_0:X_1:X_2:X_3) \mapsto (X_{\sigma(0)}:X_{\sigma(1)}:X_{\sigma(2)}:X_{\sigma(3)}),
	$$
	and if $\phi\in (\Z/2\Z)^2$, then
	$$
	\phi: (X_0:X_1:X_2:X_3) \mapsto ((-1)^{r_0}X_0:(-1)^{r_1}X_1:(-1)^{r_2}X_2:(-1)^{r_3}X_3) 
	$$
	such that $\sum_{i=1}^4 r_i\equiv 0\ (mod\ 2)$.
	
	The $H$-action on $X_{\lambda}$ gives rise to an $H$-action on $H^2_{et}(\overline{X}_{\lambda}, \overline{\Q}_{\ell}(1))=H^2_{et}(\overline{X}_{\lambda}, \Z_{\ell}(1))\otimes \overline{Q}_{\ell}$ that commutes with the $G_F$-action. So according to Proposition \ref{Canonical dec}, each factor of the canonical decomposition of $H^2_{et}(\overline{X}_{\lambda}, \overline{\Q}_{\ell}(1))$ as an $H$-representation is also a sub $G_F$-representation. By fixing a subgroup $H'<H$ and considering every canonical factor as an $H'$-representation, we can further decompose the factor canonically. As a result of the decomposition and Poincar\'e duality, we give the proof to the equation $(\ref{Eqn-main})$ and part $(a)$ of Theorem \ref{main Thm}. Then using the fact that the $G_F$-action is commutative with the  $H$-action, we prove $(c)$ and part of $(b)$ in Theorem \ref{main Thm}.

	\begin{pf}[Proof of (\ref{Eqn-main}) and (a) of Theorem \ref{main Thm}]
	In order to realize the canonical decomposition described above, we need to compute the character of $H^2_{et}(\overline{X}_{\lambda}, \overline{\Q}_{\ell}(1))$ as an $H$-representation, and find the character table of all irreducible representations of $H$. For the former, let $W$ be the $1$-dimensional subspace in $H^2_{et}(\overline{X}_{\lambda}, \overline{\Q}_{\ell}(1))$ generated by the hyperplane section, and $H^2_{pr}(X_{\lambda})$ the orthogonal complement to $W$. Let $\chi_{pr}$ the character of the sub $H$-action on $H^2_{pr}(X_{\lambda})$. Then we use the following lemma.
	
	\begin{Lem}\cite[thm.~3.2.14 ]{Gabriell-thesis} or \cite[cor.~2.5]{Gabriel-repn}. 
		If $X$ is a smooth degree $d$ projective hypersurface of dimension $n$, and $X$ is stable under the action of a projective transformation $\sigma$ of finite order, then 
		$$
		\chi_{pr}(\sigma)=\frac{(-1)^n}{d}\sum\limits_{\alpha^d=1}(1-d)^{m_{\alpha}(\sigma)}
		$$
		where $m_{\alpha}(\sigma)$ is the multiplicity of $\alpha$ as an eigenvalue of the linear representative of $\sigma$
		which leaves invariant a defining polynomial for $X$.
	\end{Lem}
	
	\begin{Cor}
		Let $\sigma\in S_{n+2}$ be such that $\sigma$ acts on $X$ by permuting the variables in its defining polynomial. Write $\sigma$ as a disjoint product of cycles, and define $m'_e(\sigma)$ to be the number of cycles whose length is divisible by $e$. Then 
		$$
		\chi_{pr}(\sigma)=\frac{(-1)^n}{d}\sum\limits_{e|4}\varphi(e)(1-d)^{m'_{e}(\sigma)}
		$$
		where $\varphi$ is the Euler totient function.
	\end{Cor}
	
	In our situation,  $X_{\lambda}$ is a surface with degree $4$. Thus
	$$
	\chi_{pr}(h)=\frac{1}{4}\sum\limits_{\alpha^4=1}(-3)^{m_{\alpha}(h)}.
	$$
	To compute $m_{\alpha}(h)$, consider the $4$-dimensional vector $\overline{\Q}$-space $U$ which is formally generated by basis $X_0,X_1,X_2,X_3$. Then according to the group action of $H$ on $U$, for every $h\in H$, $m_{\alpha}(h)$ can be interpreted as the dimension of the eigenspace of the eigenvalue $\alpha$ (if $\alpha$ is not an eigenvalue of $h$, then $m_{\alpha}(h)=0$). 
	
	\begin{eg}
		If $h=(123)=(123)(4)$, then we know that $m'_1(h)=2$, $m'_3(h)=1$, and $m'_2(h)=m_4(h)=0$. Thus $\chi_{pr}(h)=\frac{1}{4}(\varphi(1)(-3)^2+\varphi(2)(-3)^0+\varphi(4)(-3)^0)=\frac{1}{4}(9+1+2)=3$.
	\end{eg}
	\begin{eg}
		If $h(X_0:X_1:X_2:X_3)=(X_0:X_1:-X_2:-X_3)$, then $m_1(h)=2$ and $ m_{-1}(h)=2$ and $ m_{i}(h)=m_{-i}(h)=0$. Thus $\chi_{pr}(h)=\frac{1}{4}((-3)^2+(-3)^2+(-3)^0+(-3)^0)=5$.
	\end{eg}

	By similar calculations, we can write down Table \ref{Tab-1}, which is the character table of $\chi_{pr}$ (for each conjugation class, we write a representative of it, the number of elements in this class, and compute the corresponding character). And in this table, $e_i\in \Z/2\Z\subset H$ is such that $e_1(X_0:X_1:X_2:X_3)=(X_0:X_1:-X_2:-X_3)$, $e_2(X_0:X_1:X_2:X_3)=(X_0:-X_1:X_2:-X_3)$, $e_3(X_0:X_1:X_2:X_3)=(X_0:-X_1:-X_2:X_3)$. 
	\begin{table}
		\begin{center}
			\begin{tabular}{l*{10}{c}r}
				Representative                & id & $e_1$ & (12) & $e_2$(12)& (12)(34)   \\
				\hline 
				Number of elements & 1 & 3 & 12 & 12 & 3\\
				\hline
				Character $\chi_{pr}(h)$  & 21 & 5 & -7 & -3 & 5  \\
				\hline \hline
				Representative                & $e_2$(12)(34) & $e_3$(12)(34) & (123) & (1234) & $e_2$(12)(34)   \\
				\hline 
				Number of elements & 3 & 6 & 32 & 12 & 12\\
				\hline
				Character $\chi_{pr}(h)$  & 5 & 5 & 3 & -3 & -3  \\
				\hline
			\end{tabular}
			\caption{Character table of $\chi_{pr}$} \label{Tab-1}
		\end{center}
	\end{table}

	Now we compute the character table of all irreducible representations of $H$. Let $\mathcal{X}:=Hom((\Z/2\Z)^2, \overline{\Q}_{\ell}^{*})$, which has four elements $\{id, \chi_1,\chi_2,\chi_3\}$ such that $\chi_i(e_j)=\pm 1$ and equal to $1$ if and only if $i=j$. Since $H=S_4\ltimes (\Z/2\Z)^2$, we can define an $S_4$ action on  $\mathcal{X}$ by  $(s\chi)(a)=\chi(s^{-1}as)$\text{ for }$s\in S_4$ and $ \chi \in \mathcal{X}$ and $ a\in (\Z/2\Z)^2
	$. Under this action, the set $\{id, \chi_1,\chi_2,\chi_3\}$ has two orbits $[id]$ and $[\chi_1]=\{\chi_1,\chi_2,\chi_3\}$.
	By \cite[prop.~25]{Serre-LRF} the irreducible representations of $H$ bijectively correspond to the pairs $(\chi,\rho)$, where $\chi\in \mathcal{X}$ is a representative of an orbit and $\rho$ is an irreducible representation of the stable subgroup in $S_4$ to $\chi$ with respect to the above group action. For the first class $[id]$, its stable subgroup is $K_0=S_4$. As for the second class, the stable subgroup is $K_1=\{id, (12), (34), (1324), (12)(34), (1423), (13)(24), (14)(23)\}\simeq D_8$.

	With the character tables of $S_4$ and $D_8$ separately, we get a character table of $H$. To save notation, we still use $\rho$ to refer to the induced representation induced by the pair $(\chi, \rho)$.
	\begin{table}[H]
		\begin{center}
			\resizebox{\textwidth}{!}{
			\begin{tabular}{l*{10}{c}r}
				Repn/Classes                & id & $e_1$ & (12) & $e_2$(12)& (12)(34)  & $e_2$(12)(34) & $e_3$(12)(34) & (123) & (1234) & $e_2$(12)(34)   \\
				\hline 
				\# of elements & 1 & 3 & 12 & 12 & 3& 3 & 6 & 32 & 12 & 12\\
				\hline
				$\rho_1$  & 1 & 1 & 1 & 1 & 1 & 1 & 1 & 1 & 1 & 1  \\
				\hline
				$\rho_2$  & 1 & 1 & -1 & -1 & 1 & 1 & 1 & 1 & -1& -1  \\
				\hline
				$\rho_3$  & 2 & 2 & 0 & 0 & 2 & 2 & 2 & -1 & 0 & 0  \\
				\hline
				$\rho_4$  & 3 & 3 & 1 & 1 & -1 & -1 & -1 & 0 & -1 & -1  \\
				\hline
				$\rho_5$  & 3 & 3 & -1 & -1 & -1 & -1 & -1 & 0 & 1 & 1  \\
				\hline
				$\varphi_1$  & 3 & -1 & 1 & -1 & 3 & -1 & -1 & 0 & 1 & -1  \\
				\hline
				$\varphi_2$  & 3 & -1 & -1 & 1 & -1 & 3 & -1 & 0 & 1 & -1  \\
				\hline
				$\varphi_3$  & 3 & -1 & 1 & -1 & -1 & 3 & -1 & 0 & -1 & 1  \\
				\hline
				$\varphi_4$  & 3 & -1 & -1 & 1 & 3 & -1 & -1 & 0 & -1 & 1  \\
				\hline
				$\varphi_5$  & 6 & -2 & 0 & 0 & -2 & -2 & 2 & 0 & 0 & 0  \\
				\hline
			\end{tabular}
		}
		\end{center}
		\caption{Character table of $H$}
	\end{table}
	 \noindent In the above table, the $\rho_i$'s correspond to the irreducible characters of $S_4$, and the $\varphi_i$'s correspond to those of $D_{8}$. In more detail, $\rho_1$ is the trivial character, $\rho_2$ is the corresponding to the sign representation, $\rho_3$ has degree $2$, $\rho_4 $ and $\rho_5$ have degree $3$. Also, $\varphi_1$ is the trivial character,  and $\varphi_5$ is the only character of $D_8$ which has degree $2$.

	Comparing the character table of $\chi_{pr}$ with the one of $H$, we obtain the canonical decomposition:
	\begin{equation}
	H^2_{pr}(X_{\lambda})=V_{\rho_2}\bigoplus V_{\rho_3}\bigoplus V_{\rho_5}\bigoplus V_{\varphi_2}\bigoplus  V_{\varphi_4}\bigoplus  V_{\varphi_5},\label{decomposition}
	\end{equation}
	
	\noindent where $V_{\chi}$ stands for the subrepresentation corresponding to the irreducible representation $\chi$. Specifically, all but one of the $V_{\chi}$'s appearing in (\ref{decomposition}) are irreducible $H$-representations, except that $V_{\rho_2}\simeq \rho_2^{\oplus 4}$.
	
	\begin{Rmk}\label{rational dec}
		In our case, the canonical decomposition is over $\Q$. This is because all of the character values in the above tables are integers \cite[Chap.~12, Prop.~33]{Serre-LRF}. This implies that the above decomposition is also a decomposition of $H^2(\overline{X}_{\lambda}, \Q)$. Thus the decomposition induces a canonical decomposition of $H^2_{et}(\overline{X}_{\lambda}, \overline{\Q}_{\ell}(1))$.
	\end{Rmk}
	
	In fact, we can decompose every $V\neq V_{\rho_2}$ listed above. To do this, we view $V$ as a $H'$-representation for a properly chosen subgroup $H'$ of $H$, and apply the same idea as above. As an example, we show how to decompose $V_{\rho_5}$. Let $H'=S_3$. Using the character tables of $S_3$ we obtain
	$$
	V_{\rho_5}=V^{(S_3)}_2\bigoplus V^{(S_3)}_3
	$$
	Here we write $V^{(S_3)}$ to indicate that it is a $S_3$-representation. We find that $\dim V^{(S_3)}_2=1$ and $\dim V^{(S_3)}_3=2$.
	Next we choose $H'=A_3\simeq C_3$, which is a subgroup of $S_3$. As $A_3$ representation
	$$
	V^{(S_3)}_3=V^{(A_3)}_2 \bigoplus V^{(A_3)}_3
	$$
	where $V^{(A_3)}_2,V^{(A_3)}_3$ are the two nontrivial 1-dimensional  $A_3$-representations. Thus $V_{\rho_5}$ is a direct sum of three 1-dimensional $G_F$-representations.

	Similarly, we also decompose $V_{\rho_3}, V_{\varphi_2}, V_{\varphi_4}, V_{\varphi_5}$ into direct sums of 1-dimensional $G_F$-representations. More specifically:
	\begin{enumerate}
		\item $V_{\rho_3}$ is the direct sum of two distinct $1$-dimensional $A_4$ representations.
		\item $V_{\varphi_3}$ is the direct sum of three distinct $1$-dimensional $A_3$ representations.
		\item $V_{\varphi_4}$ is the direct sum of three  distinct $1$-dimensional $A_4$ representations.
		\item $V_{\varphi_5}$ is the direct sum of three distinct $1$-dimensional $A_3$ representations and a $1$-dimensional $S_3$ representation and two distinct $1$-dimensional $A_3$ representations. 
	\end{enumerate}
	
	To summarize, except for $V_{\rho_2}$, each direct summand of $(\ref{decomposition})$ is a direct sum of $1$-dimensional $G_F$-submodules. Thus we have decomposed $H^2_{et}(X_{\lambda},\overline{\Q}_{\ell}(1))$ into a direct sum of eighteen $1$-dimensional $G_F$-representations and a $4$-dimensional representation $V_{\rho_2}$. 
	
		To obtain the equation $(\ref{Eqn-main})$, we use the equation $(\ref{decomposition})$  and the summary above. If we denote those $1$-dimensional $G_F$-representations in the summary by $W_1,\cdots , W_{18}$, then we have
		$$
		H^2_{et}(X_{\lambda},\overline{\Q}_{\ell}(1))=(\bigoplus\limits_{i=1}^{18}W_i)\bigoplus V_{\rho_2}.
		$$ 
		Since for nonsingular $X_{\lambda}$, we have $rank(T_{\lambda})=3$, it suffices to show that $T_{\lambda}\subset V_{\rho_2}$. To prove this, recall that by Remark \ref{rational dec}, the canonical decompositions of $H^2(\overline{X}_{\lambda}, \C)$ and $H^2_{et}(\overline{X}_{\lambda}, \overline{\Q}_{\ell}(1))$ are compatible in the sense that they are both induced by that of $H^2(\overline{X}_{\lambda}, \Z)$. So it suffices to verify this in the complex cohomology. Take $N:=A_4\ltimes (\Z/2\Z)^2$ to be a subgroup of $H$ which consists of symplectic automorphisms (i.e. the automorphisms which preserve elements in $H^{0,2}(X_{\lambda},\C)$ \cite[page.~188, 3.3.2]{Bini-Garbagnati}). Then $N$ also preserves all elements in $H^{2,0}(X_{\lambda},\C)$ according to Poincare duality. In addition, the only nontrivial element in the quotient $H/N$ acts by multiplying $-1$ on both  $H^{2,0}(X_{\lambda},\C)$ and $H^{0,2}(X_{\lambda},\C)$. Thus we know that both $H^{2,0}(X_{\lambda},\C)$ and $H^{0,2}(X_{\lambda},\C)$ are two isomorphic 1-dimensional representations of $H$ over $\C$. Now in the canonical decomposition of $H^2_{pr}(X_{\lambda})$ with the exception of  $V_{\rho_2}$, all other canonical components are irreducible $H$-representations of dimension $>1$; thus ir follows that $H^{2,0}(X_{\lambda},\C)$ and $H^{0, 2}(X_{\lambda},\C)$ have to be contained in $V_{\rho_2}$.  Moreover, since the $3$-dimensional space $T_{\lambda}\otimes \C$ is stable under the $H$-action and contains $H^{0,2}(X_{\lambda},\C)\oplus H^{2,0}(X_{\lambda},\C)$, hence we know $T_{\lambda}$ is contained in $V_{\rho_2}$. As result, we know $T_{\rho_2}=T_{\lambda}\oplus W_{19}$ for some $1$-dimensional $G_F$-representation $W_{19}$, and obtain the following decomposition of $G_F$-representation.
		$$H^2_{et}(\overline{X}_{\lambda},\overline{\Q}_{\ell}(1))=T_{\lambda}\oplus W_1\oplus\cdots \oplus W_{19} . $$
		Next, we claim that every $W_{i}$ ($i=1,\cdots ,19$) corresponds to a quadratic character. In fact, this follows from Poincar\'e duality. Choose $v\in W_i$, then for any $g\in G_{\Q}$, we have $g(v)=\alpha v$ for some $\alpha\in \C$. Recall that the cup product is compatible with Galois action and $H^4_{et}(X,\Z_{\ell}(2))$ has trivial Galois action; if we choose $u=v$, then 
		$$
		v \cup v=g(v\cup v)=g(v)\cup g(v)=\alpha^2 v\cup v\Rightarrow \alpha^2=1
		$$
		provided that $v\cup v\neq 0$. To see why $v\cup v\neq 0$, recall that the Hodge Index Theorem (proposition \ref{ind thm} and corollary \ref{signature}) tells us that the intersection form on $NS(\overline{X}_{\lambda})$ is $(1,\rho(X)-1)$. Since the hyperplane section of $X$ has positive self-intersection number, the bilinear form restricted to $NS(\overline{X}_{\lambda})\cap H_{pr}^2(X_{\lambda})$ is negative definite. Since each $W_i$ is either equal to the space spanned by the hyperplane section or contained in $NS(\overline{X}_{\lambda})\cap H_{pr}^2(X_{\lambda})$, we know if $v\in W_i$, then $v\cup v\neq 0$. As the result of above argument, we have the following decomposition
		$$
		NS(\overline{X}_{\lambda})=id \bigoplus\limits_{i=1}^r {a_i \lrd \bullet}^{\oplus n_i}
		$$
		such that $\sum\limits_{i=1}^{r}n_i=18$, the trivial character comes from the hyperplane section, and the integers $a_i$ are distinct. This proves the equation $(\ref{main Thm})$, and $(a)$ follows immediately. 	
		\end{pf}	
	
		\begin{pf}[proof of (c) and  part of (b)]
		If we assume $U_i:= {a_i \lrd \bullet}^{\oplus n_i}$, then $H^2_{pr}(\overline{X}_{\lambda})=(\bigoplus\limits_{i=1}^r U_i)\otimes \overline{\Q}_{\lambda}$. Choose an element $v\in U_1$. For arbitrary $h\in H$, we can write $h(v)$ as a sum of its components in $U_i$: $h(v)=a_1u_1+a_2u_2+\cdots a_ku_k$ $(u_i\in U_i)$. Then, we compute $gh(v)$ for an element $g\in G$. On one side, we know it is equal to $g(a_1u_1+a_2u_2+\cdots a_ku_k)=\alpha_g^{(1)}a_1u_1+\cdots+\alpha_g^{(k)}a_ku_k$, where $\alpha_g^{(i)}:=\frac{g(u_i)}{u_i}$; on the other side, since the action of $g$ and $h$ commute, we know that $gh(v)=hg(v)=h(\alpha_g^{(1)}u)=\alpha_g^{(1)}a_1u_1+\cdots+\alpha_g^{(1)}a_ku_k$. Since different $u_i$ come from different $U_i$, when $i\neq 1$, there is some $g\in G_{F}$ such that $\alpha_g^{(1)}\neq \alpha_g^{(i)}$. This forces $a_i=0$ for all $i\neq 1$ and it follows immediately that $h(v)\in U_1$ for all $h$. Thus $U_1$ is stable under $H$-action, and so are all other $U_i$ by the same method. Recall that all of the components except $V_{\rho_2}$ in the canonical decomposition $(\ref{decomposition})$ are irreducible $H$-representations. Thus each $U_i$ contains at least one of the irreducible $H$-representations. Moreover, since (except for $V_{\rho_2}$ and $V_{\rho_3}$) every other direct summand in $(\ref{decomposition})$ has dimension divisible by $3$, and $T_{\lambda}\subset V_{\rho_2}$ has rank $3$, we see that either all $n_i$ are divisible by $3$ or exactly two of them are not divisible by $3$. This proves part of $(b)$. Finally, since $V_{\varphi_5}$ has dimension six, we know at least one of the $n_i$ is no less than $6$, and this proves $(c)$.		
	\end{pf}

	\section{The quotient family $Y_{\lambda}$}
	In this section, we finish the proof Theorem \ref{main Thm} by proving $(b)$ and $(d)$ in this theorem. In order to do this, we consider another family of $K3$ surfaces $\{Y_{\lambda}\}$ defined below. For $\lambda\in F$ such that $\lambda^4\neq 1$ and $\lambda\neq 0$, by studying the transcendental lattices of $X_{\lambda}$ and $Y_{\lambda}$, we show that they are isomorphic as $G_F$-representations. As a consequence, we finish the proof to Theorem \ref{main Thm} and also prove Corollary \ref{Cor}.
	
	For each $\lambda$, fix a primitive fourth root of unity $\xi$ and consider the group consisting of actions 
	$$
	\resizebox{1.0 \textwidth}{!} 
	{
	$A=\{\alpha:\alpha(X_0:X_1:X_2:X_3)=(\xi^{r_0}X_0:\xi^{r_1}X_1:\xi^{r_2}X_2:\xi^{r_3}X_3) \text{ such that }\sum r_i\equiv 0\ (mod\ 4)\}$.
	}
	$$
	One has $A\simeq (\Z/4\Z) ^{2}$. Denote by $M_{\lambda}$ the quotient of $X_{\lambda}$ modulo $A$. Then $M_{\lambda}$ is a variety defined by
	$$
	(Y_0+Y_1+Y_2+Y_3)^4-(4\lambda)^4Y_0Y_1Y_2Y_3=0,
	$$ 
	and the quotient map from $X_{\lambda}$ to $M_{\lambda}$ is $\phi:(X_0:X_1:X_2:X_3)\mapsto (X_0^4:X_1^4:X_2^4:X_3^4)$. After resolving the six  $A_3$ singularities of $M_{\lambda}$, we obtain a smooth surface $Y_{\lambda}$.
	
	\begin{Rmk}\label{blowupRmk}
		In fact all the six singularities of $M_{\lambda}$ are $F$-rational. And $Y_{\lambda}$ has eighteen $F$-rational exceptional curves which come from blowing-up those singular points. 
	\end{Rmk}
	
	\begin{Rmk}
		For the same $\lambda$ such that $\lambda \neq 0$ and $\lambda^4\neq 1$, the pair $(X_{\lambda}, Y_{\lambda})$ is a mirror symmetry pair (c.f.\cite{Wan-Mirror-Sym}).
	\end{Rmk}
	
	\begin{Lem}\cite{Elkies-Schutts-K3} or \cite[3.2.2]{Bini-Garbagnati}. For all $\lambda\neq 0$ and $\lambda^4\neq 1$ $Y_{\lambda}$ is also a $K3$ surface and $\rho(\overline{Y}_{\lambda})=\rho(\overline{X}_{\lambda})$.
	\end{Lem}
	
	As an immediate consequence of the construction of $Y_{\lambda}$ and Remark \ref{blowupRmk}, we have
	\begin{Cor}
		If $Y_{\lambda}$ is smooth and has Picard rank $19$, then the $G_F$-action on $NS(Y_{\lambda})$ is trivial.
	\end{Cor}
	\begin{pf}
		In fact, the eighteen lines in Remark \ref{blowupRmk} are algebraically independent since they all have self-intersection number $-2$. Moreover, those eighteen lines plus the hyperplane section become a basis of $Y_{\lambda}$ when the Picard number of $Y_{\lambda}$ is $19$. Then the conclusion follows immediately from the fact that the eighteen lines are defined over $F$.
	\end{pf}
	
	Now we consider the relation between $T(X_{\lambda})$ and $T(Y_{\lambda})$. For this, we first recall the double cover example in \cite[$\S$.~3]{Morrison-K3-SI}. Let $\iota$ be an involution of a smooth surface $X$ with isolated fixed points $Q_1, \cdots , Q_k$. Let $X\overset{\phi}{\to }X/\left< \iota \right>=:M$ be the quotient map. Denote by $\pi:Y\to M$ the minimal resolution. Also, take $\tilde{\pi}:Z\to X$ to be the blow-up of $X$ at points $Q_1, \cdots , Q_k$. Take $\tilde{\iota}:Z\to Z$ to be the involution induced by $\iota$ and $\phi: Z\to Z/\left< \tilde{\iota} \right>=Y$ to be the corresponding quotient morphism. Then the following diagram commutes.
	\begin{diagram}
		Z & \rTo^{\tilde{\pi}} & X\\
		\dTo^{\psi} & \ldDotsto^{\gamma} & \dTo^{\phi}\\
		Y & \rTo^{\pi} & M
	\end{diagram}
	
	\begin{Def}
		With the notations above, there is a degree $2$ rational map $\gamma: X\to Y$ defined by
		$$ 
		X-\{Q_1,\cdots, Q_k\}\overset{\tilde{\pi}^{-1}}{\to} Z\overset{\psi}{\to} Y.
		$$
		We call it the rational double cover map induced by $\iota$.
	\end{Def}
	
	Take $P_i=\phi(Q_i)$ and $E_i=\tilde{\pi}^{-1}(Q_i)$ and $ C_i=\pi^{-1}(P_i)$. Let $H_Z$ be the orthogonal complement of $\{E_i\}$ in $H^2(Z,\Z)$. And let $H_Y$ be the orthogonal complement of $\{C_i\}$ in $H^2(Y,\Z)$, then
	\begin{Lem}\cite[$\S$3]{Morrison-K3-SI} or \cite[appendix A, $\S$1]{GTM52} \label{Gysin}
		With the notations above, there exist natural maps
		$$
		\psi^*: H_Y\to H_Z\simeq H^2(X,\Z)
		$$
		$$
		\psi_*: H^2(X,\Z)\simeq H_Z\to H_Y \subset H^2(Y,\Z)
		$$
		such that 
		$$
		\psi_* \psi^*(y)=2y; \ \psi^*\psi_*(x)=x+\iota^*(x); \ y_1\cup y_2=\frac{1}{2} (\psi^*y_1\cup \psi^*y_2)
		$$
		Moreover, 
		$$
		\psi_*(NS(X))\subset NS(Y), \text{ and } \psi^*(K_Y)=K_X
		$$
		where $K_X$ and $K_Y$ are the canonical divisors of $X$ and $Y$, respectively.
	\end{Lem}
	
	\begin{Rmk}
		In fact, the map $\psi_*$ in lemma \ref{Gysin} is the adjoint of the pull-back map $\psi^*$. For more details of $\psi_*$, see \cite[lecture 7]{Lewis-Survey-Hodge} and \cite[chap.~\RNum{6} $\S$.~5]{MilEC}.
	\end{Rmk}
	
	Note that we have a filtration of groups $A=(\Z/4\Z)^2=H_0\supset H_1\supset H_2\supset H_3\supset H_4=0$ such that $H_{i}/H_{i+1}\simeq \Z/2\Z$. Moreover, the fixed points of $X_{\lambda}$ by $A$ are isolated. Thus starting with $X_{\lambda}$ for a fixed $\lambda$, the process of constructing $Y_{\lambda}$ is a composition of four rational double cover maps induced by the quotients $A_i/A_{i+1}$ for $i=0,\cdots, 3$. As an analogy of the double cover case, we have
	
	\begin{Prop}
		For a fixed smooth $X_{\lambda}$, there exists a smooth surface $Z_{\lambda}$, birational surjection $\tilde{\pi}: Z_{\lambda}\to X_{\lambda}$. The $(\Z/4\Z)^2$ group action on $X_{\lambda}$ induces a $(\Z/4\Z)^2$ group action on $Z_{\lambda}$. Modulo $(\Z/4\Z)^2$ there is a degree $16$ quotient morphism $\psi: Z_{\lambda}\to Y_{\lambda}$, such that the following diagram commutes
		\begin{diagram}
			Z_{\lambda} & \rTo^{\tilde{\pi}} & X_{\lambda}\\
			\dTo^{\psi} &  & \dTo^{\phi}\\
			Y_{\lambda} & \rTo^{\pi} & M_{\lambda}
		\end{diagram}
	\end{Prop}
	Similarly, if we denote by $\{P_i\}$ the set of singular points on $M_{\lambda}$, $C_i=\pi^{-1}(P_i)$, and let $H_{Y_{\lambda}}$ be the orthogonal complement of $\{C_i\}$ in $H^2_{et}(\overline{Y}_{\lambda},\Q_{\ell}(1))$, then we have a corollary to Lemma \ref{Gysin} in \'{e}tale cohomology.
	
	\begin{Cor}
		There exist natural maps
		$$
		\psi^*: H_{Y_{\lambda}}\to  H^2_{et}(\overline{X}{_{\lambda}},\Q_{\ell}(1)) \text{ and }
		$$
		$$
		\psi_*: H^2_{et}(\overline{X}{_{\lambda}},\Q_{\ell}(1))\to H_{Y{_{\lambda}}}
		$$
		such that 
		$$
		\psi_* \psi^*(y)=16y \text{ and } \psi^*\psi_*(x)=\sum\limits_{\alpha\in A}\alpha^*(x)  \text{ and }  y_1\cup y_2=\frac{1}{16} (\psi^*y_1\cup \psi^*y_2)
		$$
		Moreover, $\psi^*(T_{Y_{\lambda}}\otimes \Q_{\ell})=T_{X_\lambda}\otimes \Q_{\ell}$.
	\end{Cor}
	\begin{pf}
		Everything except the last statement follows from Lemma \ref{Gysin}. Note that $\psi_*\psi^*=16$ implies that $\psi^*$ is injective. Also note the fact that $\psi_*$ is the adjoint of $\psi^*$ with respect to the cup product, so $\psi^*(T_{Y_{\lambda}})$ is contained in $T_{X_{\lambda}}\otimes \Q_{\ell}$. Thus the last statement is true due to the fact that $T_{X_{\lambda}}$ and $T_{Y_{\lambda}}$ have the same rank.
	\end{pf}
	
	\begin{Cor}\label{T have the same repn}
		The transcendental lattices $T_{X_\lambda}$ and $T_{Y_{\lambda}}$ have the same $G_F$-action. 
	\end{Cor}  
	
	\begin{pf}
		It is sufficient to show that $\psi^*$ commutes with the $G_F$-action. But this follows from the fact both $\pi_*$ and $\phi^*$ do.
	\end{pf}
	
	Now we can prove corollary \ref{Cor}.  In fact, we can prove a stronger result. Let $\mathcal{O}_K$ be the algebraic integer ring of $K$ and let $\mathfrak{p}$ be a prime ideal $\mathcal{O}_K$ such that both $X_{\lambda} $ and $Y_{\lambda}$ have good reduction at $\mathfrak{p}$. Assume that the norm $|\mathcal{O}_K/\mathfrak{p}|$ is a power of prime number $p$ and $q=|\mathcal{O}_K/\mathfrak{p}|^k$. 
	
	\begin{Lem}\label{Lem}
		Using the notations above and taking the assumption that $\rho(\overline{X}_{\lambda})=19$, we have $\#X_{\lambda}(\F_q)\equiv \#Y_{\lambda}(\F_q)\ (mod\ 3q)$ when $p\neq 2$ or $3$.
	\end{Lem}
	\begin{pf}
		Firstly, note that 
		
		\begin{equation}
			\resizebox{1.0\textwidth}{!}{%
			  $\begin{aligned}
			  \#(X_{\lambda}(\F_q))-\#(Y_{\lambda}(\F_q)) &=tr(Frob_q|H^2_{et}({\overline{X}}_{\lambda},\overline{Q}_{\ell}(1)))-tr(Frob_q|H^2_{et}({\overline{Y}}_{\lambda},\overline{Q}_{\ell}(1)))\\
			  		& =tr(Frob_q|NS({\overline{X}}_{\lambda})\otimes \Q_{\ell})-tr(Frob_q|NS({\overline{Y}}_{\lambda})\otimes \Q_{\ell}).	
			  \end{aligned}$%
			}
			\end{equation}
		
		Hence $\#X_{\lambda}(\F_q)\equiv \#Y_{\lambda}(\F_q)\ (mod\ q)$ follows immediately from this equation, equation $(\ref{Eqn-main})$, and the fact that $NS(\overline{Y}_{\lambda})$ has trivial $G_F$-action. So it is reduced to show  $\#X_{\lambda}(\F_q)\equiv \#Y_{\lambda}(\F_q)\ (mod\ 3)$.
		
		Secondly, notice that $Y_{\lambda}$ comes from blowing-up $M_{\lambda}$ at its six singularities (cf. Remark \ref{blowupRmk}). Thus this problems reduces to showing that $\#X_{\lambda}(\F_q)\equiv \#M_{\lambda}(\F_q)\ (mod\ 3)$. Moreover, by the relationship of a projective variety and its affine cone, this problem is further reduced to showing that their affine cones $A(X_{\lambda})$,$A(M_{\lambda})$ have the same number of points $mod\ 3$.  Both of $A(X_{\lambda})$ and $A(M_{\lambda})$ admit action from $S_4$. For any point $P=(x_0,x_1,x_2,x_3)\in A(X_{\lambda})(\F_q)$ or $P=(y_0,y_1,y_2,y_3)\in A(M_{\lambda})(\F_q)$, we use $[P]$ to denote $S_4$ orbit of $P$ and denote by $\#[P]$ the number of elements in this orbit. 
		
		Thirdly, we count $A(X_{\lambda})(\F_q)\ (mod\ 3)$. To do this, note if at least three of $x_0,x_1,x_2,x_3$ are distinct, then $\#[P]=\#[(x_0,x_1,x_2,x_3)]\equiv 0 \ (mod\ 3)$. Then if $P=(a,a,b,b)$  and $a\neq b$ with $ab\neq 0$, we can see that $\#[P]=\binom{4}{2}=6\equiv 0\ (mod\ 3)$. If $P=(a,b,b,b)$ with $a\neq b$ and $ab\neq 0$, the point $(a,b)$ is on the curve $C_{X_{\lambda}}:u^4+3v^4-4\lambda uv^3=0$. At last, if $P=(a,a,a,a)$ then we have $4a^4-4\lambda a^4=0$. Since we require $\lambda^4\neq 1$, $a$ has to be $0$, and $P=(0,0,0,0)$ is in $A(X_\lambda)$. 
				
		Next, we apply the same idea to $A(M_{\lambda})(\F_q)$. So we note that when $Q$ has at least three distinct coordinates, or when $Q=(a,a,b,b)$ with $a\neq b$ and $ab\neq 0$, we always have $\#[Q]\equiv 0\ (mod\ 3)$. When $Q=(a,b,b,b)$ with $ab\neq 0$ and $a\neq b$, then the point $(a,b)$ is on the curve  $C_{Y_{\lambda}}:(s+3t)^4-(4\lambda)^4st^3=0$. If $Q=(a,a,a,a)$, by a similar arguments, we see that $a=0$, and $Q=(0,0,0,0)\in A(Y_{\lambda})(\F_q)$. 
		
		Finally, we are reduced to comparing the number of $\F_q$ points on two curves $C_{X_{\lambda}}$, and $C_{Y_{\lambda}}$. Since we require $uv, st\neq 0$, it remains to compare $u^4+3-4\lambda u=0$ and $(s+3)^4-(4\lambda)^4s=0$. Then one notices that if $\alpha$ is a root of $u^4+3-4\lambda u=0$, then $\alpha^4$ is a root of $(s+3)^4-(4\lambda)^4s=0$ resulting from the quotient map from $X_{\lambda}$ to $M_{\lambda}$. Moreover, when $\lambda\neq 0$ (notice the case $\lambda=0$ is singular),  the numbers $-\alpha,\pm\sqrt{-1}\alpha$ are not roots. On the other hand, if $\beta$ is a root of $(s+3)^4-(4\lambda)^4s=0$, then $\frac{\beta+3}{4\lambda}$ is a root of $u^4+3-4\lambda u=0$. Thus there is a bijection between the $\F_q$ points on $C_X$ and the $\F_q$ points on $C_Y$. So conclude that $\#A(X_{\lambda})(\F_q)\equiv \#A(M_{\lambda})(\F_q)\ (mod\ 3)$ and complete the proof.		
	\end{pf}
	
	Now Corollary \ref{Cor} follows from the above lemma. In fact, with the help of the above Lemma, we can also finish the proof of our Theorem \ref{main Thm}.
	\begin{Cor}
		In equation $(\ref{Eqn-main})$, every $n_i$ is divisible by $3$. In particular, the index number $r\leq 5$.
	\end{Cor}
	\begin{pf}
		If not, according to last part of the proof in the previous section, exactly two of them are not divisible by $3$. However, it means that for some $p\neq 3$, the equation $\#X_{\lambda}(\F_p)\equiv \#Y_{\lambda}(\F_p)\ (mod\ 3)$ does not hold, contradiction.
	\end{pf}

	\section{Degree two quotients of the Dwork family}
	In the second part, we use the geometric properties of degree $2$ del Pezzo surfaces to find an explicit set of generators of the Dwork family $X_{\lambda}$. We then determine the Galois representation by direct computation. Before stating our result in the next section, we review necessary background in this section. 
	
	A del Pezzo surface is a non-singular projective algebraic surface with ample anticanonical divisor. The degree of a del Pezzo surface is defined to be the self-intersection number of its canonical divisor. Geometrically, a del Pezzo surface is isomorphic to $\mathbb{P}^2$, to $\mathbb{P}^1\times \mathbb{P}^1$, or to the blowing-up of $\mathbb{P}^2$ up to eight points in general position (\cite[$\S$ 1]{Yuri-deg-2-del-pezzo}, \cite[Thm.~8.15]{Dolgachev-Classical-Alg-Geo}). In the last situation, the degree $d$ equals $9$ minus the number of points in the blow-up. It follows that the Picard group of this surface is isomorphic to $\Z^{10-d}$ \cite[chap.~\RNum{5}, Prop.~3.2]{GTM52}. From now on, we focus on the degree two del Pezzo surfaces constructed as degree two quotients of our Dwork surfaces. The construction of these quotients follows from Kloosterman's work \cite{Kloosterman}.
	
	Fix $I$ to be a square root of $-1$, we consider the involution of the form 
	$$\tau_{i,j}^{r}: X_{i}\mapsto I^rX_j, X_{j}\mapsto (-I)^rX_i, \text{ and } X_k\mapsto X_k \text{ if } k\neq i,j$$
	 and denote by $S_{\lambda}^{(i,j;r)}$ the quotient surface $X_{\lambda}/\langle \tau^{r}_{i,j}\rangle$.
	
	\begin{eg}
	By the above construction, we have 
	\begin{align*}	
		\tau^{(4)}_{0,1}: (X_0: X_1:X_2 :X_3) & =(X_1: X_0: X_2: X_3),	\\
		\tau^{(2)}_{0,1}: (X_0: X_1:X_2 :X_3) & =(-X_1: -X_0: X_2: X_3),	\\
		\tau^{(1)}_{0,1}: (X_0: X_1:X_2 :X_3) & =(IX_1: -IX_0: X_2: X_3).
	\end{align*}
	If we let $v=X_0X_1$, and let 
	$$
	u=
	\begin{cases*}
	X_0+X_1 & \text{if }r=4,\\
	X_0-X_1 & \text{if }r=2,\\
	X_0+IX_1 & \text{if }r=1.\\
	\end{cases*}
	$$ 
	then we can find the explicit defining equations for the corresponding quotient surfaces
	\begin{align*}
	S_{\lambda}^{(0,1;4)} &: u^4-4u^2v+2v^2-4\lambda v X_2X_3+X_2^4+X_3^4=0,\\
	S_{\lambda}^{(0,1;2)} &: u^4+4u^2v+2v^2-4\lambda v X_2X_3+X_2^4+X_3^4=0,\\
	S_{\lambda}^{(0,1;1)} &: u^4-4Iu^2v-2v^2-4\lambda v X_2X_3+X_2^4+X_3^4=0.\\
	\end{align*}
	The other $S^{(i,j;r)}$ are obtained similarly. 
	\end{eg}
	
	For each $(i,j;r)$, the surface $S^{(i,j;r)}_{\lambda}$ is a degree two del Pezzo surface in the weighted projective space $\mathbb{P}(1,1,1,2)$ (\cite[proof of Prop.~3.1 and Appindex A]{Kloosterman}). Note that $S_{\lambda}^{(i,j;r)}$ is also a double covering above $\mathbb{P}^2$ which is branched over a quartic curve, thus $S^{(i,j;r)}_{\lambda}$ is smooth if and only if the branched curve is smooth, i.e. if and only if $\lambda^4\neq 1$. Moreover, when the quartic curve is smooth, it has exactly $28$ bitangents, which give rise to all the $56$ lines after pulling back to $S^{(i,j;r)}_{\lambda}$. The $56$ straight lines generate the N\'eron-Severi group of $S^{(i,j;r)}_{\lambda}$. Indeed, any seven pairwise disjoint lines plus the anticanonical divisor will be a basis of the N\'eron-Severi group (recall that a degree $2$ del Pezzo surface has Picard rank eight). For more information about degree two del Pezzo surfaces, one is referred to \cite[$\S$ 8.2, 8.7]{Dolgachev-Classical-Alg-Geo}.

	\begin{eg}\label{eg-2}
	Consider $w=v-u^2-\lambda X_2X_3$. Then we can rewrite the equation of $S^{(0,1;4)}_{\lambda}$ as 
	$$
	w^2=u^4-X_2^4-X_3^4+4\lambda u^2X_2X_3+2\lambda^2X_2^2X_3^2.
	$$
	Then a straight line $\alpha u+\beta X_2+\gamma X_3=0$ is a bitanget of the branch quartic curve $u^4-X_2^4-X_3^4+4\lambda u^2X_2X_3+2\lambda^2X_2^2X_3^2=0$ if and only if their intersection can be represented as the square of a polynomial, or equivalently, if and only if it is one of the following $7$ types. 
	\begin{enumerate}
	\item $X_2=aX_3$, where $a$ satisfies equation $a^4+2\lambda^2a^2+1=0$. Calculation shows that $a=\delta \sqrt{-2(\lambda^2-1)}+\gamma\sqrt{-2(\lambda^2+1)}\in \Q(\sqrt{-2(\lambda^2-1)}, \sqrt{-2(\lambda^2+1)})$ with $ \delta, \gamma\in \Q^*:=\Q-\{0\}$. And in this case, the two lines   $S^{(0,1;4)}_{\lambda}$ above this bitangent are given by equations $4w=\pm \sqrt{2}(4a\lambda X_3^2+2u^2)$. Moreover, for each $\lambda$, we fix the following two lines
	\begin{enumerate}
	\item $v_1$ is the line uniquely determined by $a=\sqrt{-\lambda^2-\sqrt{\lambda^4-1}}$ and $4w=-\sqrt{2}(4a\lambda X_3^2+2u^2)$, and
	\item $v_2$ is the line determined by $a=\sqrt{-\lambda^2+\sqrt{\lambda^4-1}}$ and $ 4w= \sqrt{2}(4a\lambda X_3^2+2u^2)$.
	\end{enumerate}  
	
	\item $u=-\frac{1}{a}X_2+aX_3$, and $a^4=1$. In this case, the lines above the bitangets are given by $ 2w^2=\pm X_2^2\sqrt{2(\lambda^4-1)} $. We take 
	\begin{enumerate}
	\item $v_3$ to be the line with $a=-I$ and $2w=X_2^2\sqrt{2(\lambda^4-1)}$, and
	\item $v_4$ to be the line with $a=1$ and $2w=-X_2^2\sqrt{2(\lambda^4-1)}$.
	\end{enumerate}
	
	\item $u=aX_2-\frac{1}{a}X_3$ with $a^4=1$. In this case the lines above the bitangets are given by $ 2w^2=\pm X_3^2\sqrt{2(\lambda^4-1)} $. We take
	\begin{enumerate}
	\item $v_5$ to be the line with $a=-1$ and $2w=-X_3^2\sqrt{2(\lambda^4-1)}$, and 
	\item $v_4$ to be the line with $a=I$ and $2w=X_2^3\sqrt{2(\lambda^4-1)}$, and
	\item $v'_8$ to be the line with $a=-I$ and $2w=X_3^3\sqrt{2(\lambda^4-1)}$.
	\end{enumerate}
	
	\item $u=aX_2+aX_3$ with $a^4+\frac{4\lambda}{\lambda^2+1}a^2+1$. In this case we have $a=\delta\sqrt{2(\lambda^2+1)}+\gamma \sqrt{-2(\lambda^2+1)}$. The two lines above the bitangents are given by $4w=\pm \frac{1}{\sqrt{-(\lambda^2-1)}}F$ with $F\in \Q(\lambda, a)[u,X_2,X_3]$. We fix $v_7$ to be the line with $a=-\sqrt{\frac{-2\lambda-I(\lambda^2-1)}{\lambda^2+1}}$ such that it is disjoint with $v_1$.
	
	\item $u=aX_2-aX_3$ with $a^4-\frac{4\lambda}{\lambda^2+1}a^2+1$.
	
	\item $u=aX_2+IaX_3$ with $a^4+\frac{4\lambda}{\lambda^2-1}a^2+1$.
	\item $u=aX_2-IaX_3$ with $a^4-\frac{4\lambda}{\lambda^2-1}a^2+1$.
	\end{enumerate}
	\end{eg}
	
	Now we claim that the lines $v_i$ $(i=1,\cdots, 7)$ chosen above are pairwise disjoint and hence the set which consists of them and the anticanonical divisor $-K_S$ is a basis of the N\'eron-Severi group $NS(S_{\lambda}^{(0,1,;4)})\otimes \Q$. To see this we can use algebra computing packages, or the following proposition. Let $T$ be a an irreducible variety of dimension $m>0$, and let $\mathcal{Y}$ be a smooth over $T$ of relative dimension $n$. If we take $\mathcal{Y}_t$ to be the fiber of $t\in T$, then any $(k+m)$-cycle on $\mathcal{Y}$ or more generally any rational equivalence class $\alpha\in A_{k+m}\mathcal{Y}$ determines a family of $k-$cycle classes $\alpha_t\in A_k(\mathcal{Y}_t)$ for all $t\in T$ (\cite[$\S$.~ 10.1]{Fulton-int-theory}). 
	
	\begin{Prop}\label{Prop-int}\cite[Coro.~10.1]{Fulton-int-theory}
	Assume that $T$ is non-singular and $t\in T$ is rational over the ground field. If $\alpha\in A_{k+m}(\mathcal{Y}) $ and $\beta\in A_{l+m}(\mathcal{Y})$, then 
	$$
	\alpha_t . \beta_t=(\alpha.\beta)_t
	$$
	in $ A_{k+l-n}(Y_t)$.
	\end{Prop} 
	
	Let $T=A^1-\{\pm 1,\pm I\}$ and let $\mathcal{Y}=S^{(0,1;4)}$ be restricted on $T$. Then by Proposition \ref{Prop-int}, we know that the intersection numbers $v_i.v_j$ ($i,j=1,\cdots, 7$) are independent of the choice of $\lambda$. In particular, as $\lambda$ goes to $0$, $S^{(0,1;4)}_{\lambda}$ deforms to $S^{(0,1;4)}_0$ and the intersection numbers of curves in the former surface can be calculated by the intersection numbers of the corresponding lines in the latter.
	
	\begin{eg}\cite[$\S$.2]{Yuri-deg-2-del-pezzo}\label{eg-1}
	Let $S$ be the surface defined by the equation
	$$
		2w^2=Ax^4+By^4+Cz^4.
		$$
	Fix an eighth root of unity $\zeta$, such that $\zeta^2=I$ and fix numbers $a,b,c$ such that $(a^4:b^4:c^4)=(A:B:C)$. Consider the following types of lines in $S$ (in the followings, $\delta^4=-1$ and $ \alpha^4=\beta^4=\gamma^4=1$).
	
	\begin{align*}
	 L_{z,\delta, \pm}: & \delta ax+by=0 , \ \ \ \ \ \ \ \ \ \ \ \ \sqrt{2}w=\pm c^2z^2,\\
	 L_{x,\delta, \pm}: & \delta by+cz=0 , \ \ \ \ \ \ \ \ \ \ \ \  \sqrt{2}w=\pm a^2x^2,\\
	 L_{y,\delta, \pm}: & \delta cz+ax=0 , \ \ \ \ \ \ \ \ \ \ \ \  \sqrt{2}w=\pm b^2y^2,\\
	 L_{\alpha, \beta,\gamma}: & \alpha ax+\beta by+\gamma cz=0 ,\ \  w=\alpha \beta axby+\beta \gamma bycz+\alpha \gamma czax=0.
	\end{align*}
	As $S^{(0,1;4)}_{\lambda}$ deforms to $S^{(0,1;4)}_{0}$, our lines $v_1, \cdots, v_7, v_8'$ deform to eight lines in $S^{(0,1;4)}_{0}$, which we still denote by $v_1,\cdots, v_7$ and $v'_8$. Moreover if we take $a=1$ and $ b=c=\zeta$, we have 
	\begin{align*}
	&v_1=[L_{x,\zeta, +}] \ \ \ \ v_2=[L_{x,\zeta^3, -}]  \ \ \ \ v_3=[L_{y,\zeta, +}]  \ \ \ \ v_4=[L_{y,\zeta^3, -}] \\
	&v_5=[L_{z,\zeta, +}] \ \ \ \ v_6=[L_{z,\zeta^3, -}] \ \ \ \ v_1=[L_{I,I,I}] \ \ \ \ \ v_8'=[L_{z,\zeta^7, -}] .
	\end{align*}
	Now let $v_8=v_6+v_7+v_8'$. It is easy to verify that 
	$$
	v_i.v_j=
	\begin{cases*}
	0 & if      $i\neq j$\\
	-1 & if $i=j$ and $i<8$\\
	1 & if $i=j=8$
	\end{cases*}
	$$
	and $-K_{S}=3v_8-\sum_{i=1}^{7}v_i$. Thus our claim is true, i.e. we have found a concrete basis of the N\'eron-Severi group of $S^{(0,1;4)}_{\lambda}$.
	\end{eg}
	
	\begin{Rmk}
	When deforming from $S^{(0,1;4)}_{\lambda}$ to $S^{(0,1;4)}_{0}$, our image of $v_i$ depends on the path of $\lambda$. For example, in the defining equations of $v_i$, the value of $\sqrt{-2\lambda^2-1}$ depends on our choice of the path. But we can eliminate the ambiguity by fixing a path. Moreover, the algebraic independence of those lines is not related to the path. In fact, to compute the Galois action, we do not really need to know the exact image of $v_i$ in $S^{(0,1;4)}_{0}$.
	\end{Rmk}

	\section{Explicit computation of  Galois actions}
	In this section, we make use of the double quotients $S^{(i,j;r)}_{\lambda}$ above and the generators of their N\'eron-Severi groups to determine the Galois action on the N\'eron-Severi group of the Dwork family. First, we fix the following five quotients:
	\begin{align*}
	&S_{1,\lambda}:=S^{(0,1;4)}_{\lambda}, \ \ \ \ S_{2,\lambda}:=S^{(0,1;2)}_{\lambda}, \ \ \ \ S_{3,\lambda}:=S^{(1,2;4)}_{\lambda},\\
	&S_{4,\lambda}:=S^{(2,3;4)}_{\lambda}, \ \ \ \ S_{5,\lambda}:=S^{(0,1;1)}_{\lambda}.
	\end{align*}

	By the similar method we used in the previous section, we are able to find the generators of the N\'eron-Severi groups of each $S_{i,\lambda}$. In particular, for each $i$, we denote by $v_j^{i}$ ($j=1,\cdots ,8$) the divisor which deforms to $v_j$ listed in example \ref{eg-1}. We also use the same notations for their image in $NS(\overline{X}_{\lambda})$ through pulling back. 	
	\begin{Lem}\label{Lem-S1}
	\begin{equation*}
	\resizebox{1.0 \textwidth}{!} 
	{
		$NS(\overline{S}_{1,\lambda})\simeq id \oplus {-(\lambda^2-1) \lrd \bullet}^{\oplus 1} \bigoplus {-(\lambda^2+1) \lrd \bullet}^{\oplus 1}\bigoplus {-2(\lambda^4-1) \lrd \bullet}^{\oplus 2}\bigoplus {2(\lambda^4-1) \lrd \bullet}^{\oplus 3}$
		}		
	\end{equation*}
	\end{Lem}
	
	\begin{pf}
	Take $L=\Q(\sqrt{-1},\sqrt{2},\sqrt{\lambda^2+1},\sqrt{\lambda^2-1})$. By the construction of the lines $v_j^1$ in example \ref{eg-2}, we find that they are fixed by Galois group $Gal(\overline{Q}/L)$. Now consider the set of four algebraic numbers $\{I, \sqrt{2}, \sqrt{\lambda^2+1}, \sqrt{\lambda^2-1}\}$ and take $\sigma_{I}, \sigma_2, \sigma_{+}, \sigma_{-}$ to be the elements in $Gal(L/\Q)$ which change the signs of $I, \sqrt{2}$,$ \sqrt{\lambda^2+1}, \sqrt{\lambda^2-1}$ respectively and fix the signs of all other three. Since $Gal(L/\Q)\simeq (\Z/2\Z)^{4}$ and is generated by the four elements, we are reduced to determine the actions of the four elements. To do this, we consider the actions of the four elements on each of the divisors $v_j^1$, and then deform to $\lambda=0$ to find the linear combination of the basis $\{v_{j}^1\}$  ($j=1,\cdots, 8$) which represents the result.
	
	First we determine the action of $\sigma_{I}$. By example \ref{eg-2}, the elements $v_1^1, v_2^1$ are determined by the choice of the root of $x^4+2\lambda^2x^2+1=0$ (we denote its root by $a$) and by the signs of the equation $4w=\pm \sqrt{2}(4a\lambda X_3^2+2u^2)$. Moreover, we know that the root $a$ is of the form 
	$$
	a=\delta\sqrt{-2(\lambda^2-1)}+\gamma \sqrt{-2(\lambda^2+1)},\ \ \ (\delta, \gamma\in \Q^{*}).
	$$
	Thus $\sigma_{I}(a)=-a$ and $\sigma_{I}$ fixes the sign of $4w=\pm \sqrt{2}(4a\lambda X_3^2+2u^2)$ since $\sigma_{I}(\sqrt{2})=\sqrt{2}$. This implies that after deformation (no matter which path we choose), $[L_{x,\zeta, +}]=v_1^1$ is mapped to $[L_{x,\zeta^5,+}]=-v^1_1-v^1_7+v^1_8$. Similarly, $[L_{x,\zeta^3, +}]=v_2^1$ is mapped to $[L_{x,\zeta^7,-}]=-v^1_2-v^1_7+v^1_8$.
	
	Since $v_3^1$ is determined by $a=-I$ and equation $2w=X_2^2\sqrt{2(\lambda^4-1)}$, and $v_4^1$ is determined by $a=1$ and $2w=-X_2^2\sqrt{2(\lambda^4-1)}$, we find that $v_3^1$ is mapped to the line corresponding to $a=I$, and $2w=X_2^2\sqrt{2(\lambda^4-1)}$. So after deformation, we find that $[L_{y,\zeta, +}]=v_3^1$ is mapped to $[L_{y,\zeta^5, +}]=-v_3^1-v_7^1+v_8^1$. It is easy to see that $\sigma_{I}(v_4^1)=v_4^1$.
	
	By the same argument, we see that $v_5^1$ is fixed by $\sigma_{I}$. $[L_{y,\zeta^3,-}]=v_6^1$ is mapped to $[L_{z,\zeta^7, -}]=-v_6^1-v_7^1+v_8^1$, and $[L_{z,\zeta^7,-}]={v_8^1}'$ is mapped to $[L_{z,\zeta^3, -}]=v_6^1$.
	
	At last, we see that $v_7^1$ is determined by the choice of the root $a$ of equation $x^4+\frac{4\lambda}{\lambda^2+1}x^2+1=0$ and by the sign of $\pm \sqrt{-(\lambda^2+1)}$. We can use the previous method to find $\sigma(v_7^1)$, but it is easier to do this if we note that the Galois actions preserve the intersection numbers. Thus $v_7^1$ is mapped to the only line determined by a root of $x^4+\frac{4\lambda}{\lambda^2+1}x^2+1=0$ and by the fact that it is disjoint with all the $\sigma_{I}(v_j^1)$ for $j=1,\cdots ,6$. Thus $[L_{i,i,i}]=v_7^1$ is mapped to $[L_{i,1,1}]=-K_{S}+v_4^1+v_5^1-v_7^1=-v_1^1-v_2^1-v_3^1-v_6^1-2v_7^1+3v_8^1$. Hence we have that $v_8^1$ is mapped to $-K_S+v_4^1+v_5^1-v_7^1=-v_1^1-v_2^1-v_3^1-v_6^1-2v_7^1+3v_8^1$. So with the basis $\{v_1^1, \cdots, v_8^1\}$, we can write down the matrix $M(\sigma_{I})$ corresponding to $\sigma_{I}$	
	$$
	M(\sigma_{I})=\left( \begin{array}{cccccccc}
	-1 & 0 & 0 & 0 & 0 & 0 & -1 & -1 \\ 
	0 & -1 & 0 & 0 & 0 & 0 & -1 & -1 \\ 
	0 & 0 & -1 & 0 & 0 & 0 & -1 & -1 \\ 
	0 & 0 & 0 & 1 & 0 & 0 & 0 & 0 \\ 
	0 & 0 & 0 & 0 & 1 & 0 & 0 & 0 \\ 
	0 & 0 & 0 & 0 & 0 & -1 & -1 & -1 \\ 
	-1 & -1 & -1 & 0 & 0 & -1 & -1 & -2 \\ 
	1 & 1 & 1 & 0 & 0 & 1 & 2 & 3
	\end{array} \right) 
	$$	
	By a similar method, we can write down the matrices corresponding to the remaining three actions. 	
	$$
	M(\sigma_{2})=
	\left(\begin{array}{cccccccc}
	0 & -1 & -1 & -1 & -1 & -1 & 0 & -2 \\ 
	-1 & 0 & -1 & -1 & -1 & -1 & 0 & -2 \\ 
	-1 & -1 & -2 & -1 & -1 & -1 & -1 & -3 \\ 
	-1 & -1 & -1 & -2 & -1 & -1 & -1 & -3 \\ 
	-1 & -1 & -1 & -1 & -2 & -1 & -1 & -3 \\ 
	-1 & -1 & -1 & -1 & -1 & -2 & -1 & -3 \\ 
	0 & 0 & -1 & -1 & -1 & -1 & -1 & -2 \\ 
	2 & 2 & 3 & 3 & 3 & 3 & 2 & 7
	\end{array} 
	 \right) 
	 $$
	 
	 $$
	 M(\sigma_{+})=
	 \left(
	 \begin{array}{cccccccc}
	 -1 & 0 & -1 & -1 & -1 & -1 & 0 & -2 \\ 
	 0 & -1 & -1 & -1 & -1 & -1 & 0 & -2 \\ 
	 -1 & -1 & -2 & -1 & -1 & -1 & -1 & -3 \\ 
	 -1 & -1 & -1 & -2 & -1 & -1 & -1 & -3 \\ 
	 -1 & -1 & -1 & -1 & -2 & -1 & -1 & -3 \\ 
	 -1 & -1 & -1 & -1 & -1 & -2 & -1 & -3 \\ 
	 0 & 0 & -1 & -1 & -1 & -1 & -1 & -2 \\ 
	 2 & 2 & 3 & 3 & 3 & 3 & 2 & 7
	 \end{array}  \right) 
	 $$
	 
	 $$
	 M(\sigma_{-})=
	 \left( 
	 \begin{array}{cccccccc}
	 -1 & -2 & -1 & -1 & -1 & -1 & -1 & -3 \\ 
	 -2 & -1 & -1 & -1 & -1 & -1 & -1 & -3 \\ 
	 -1 & -1 & -2 & -1 & -1 & -1 & -1 & -3 \\ 
	 -1 & -1 & -1 & -2 & -1 & -1 & -1 & -3 \\ 
	 -1 & -1 & -1 & -1 & -2 & -1 & -1 & -3 \\ 
	 -1 & -1 & -1 & -1 & -1 & -2 & -1 & -3 \\ 
	 -1 & -1 & -1 & -1 & -1 & -1 & -2 & -3 \\ 
	 3 & 3 & 3 & 3 & 3 & 3 & 3 & 8
	 \end{array} \right) 
	 $$
	 
	 Finally, the conclusion of this lemma follows by direct calculation. 	
	\end{pf}
	
	\begin{Cor}
	The N\'eron-Severi groups of $S_{2,\lambda},S_{3,\lambda}$ and $S_{4,\lambda}$ have the same decomposition as that of $S_{2,\lambda}$.
	\end{Cor}
	\begin{pf}
	This follows directly from the fact that each of the three surfaces is isomorphic to $S_{1,\lambda}$ over $F$. 
	\end{pf}
	
	Now we fix the ordered set 
	\begin{align*}
	B:=&\{e_1,e_2,e_3,e_4,e_5,e_6,e_7,e_8,e_9,e_{10},e_{11},e_{12},e_{13},e_{14},e_{15},e_{16},e_{17},e_{18},e_{19}\}\\
	=&\{v_1^1,v_2^1,v_3^1,v_4^1,v_5^1,v_6^1,v_7^1,v_8^1,v_3^2, v_4^2,v_5^2, v_6^2, v_1^3, v_2^3, v_3^3, v_4^3, v_6^3, v_1^4, v_1^5\}.
	\end{align*}
	By calculation we can see that $B$ is a basis of $NS(X_{\lambda})\otimes \Q$. Then by studying the $G_F$-action on the N\'eron-Severi groups of $S_{1,\lambda}, \cdots, S_{4,\lambda}$, we find
	
	\begin{Prop}
	As a sub-representation of $NS(\overline{X}_{\lambda})$, we have 	
	\begin{equation}
	\resizebox{1.0\textwidth}{!}{%
	  $\begin{aligned}
	  NS(\overline{S}_{1,\lambda})+NS(\overline{S}_{2,\lambda})
	   & \simeq id \oplus {-(\lambda^2-1) \lrd \bullet}^{\oplus 1}   \bigoplus {-(\lambda^2+1) \lrd \bullet}^{\oplus 1}\\
	   &\bigoplus {-2(\lambda^4-1) \lrd \bullet}^{\oplus 4}\bigoplus {2(\lambda^4-1) \lrd \bullet}^{\oplus 5}
	  \end{aligned}$%
	}
	\end{equation}

	\end{Prop}	
	As a consequence of Theorem \ref{main Thm}, (b), we obtain Theorem \ref{new main Thm}, and thus finish the proof.
	
	\begin{Rmk}
	Without using Theorem \ref{main Thm} we are still able to prove our main theorem. To do this, we need the explicit Galois actions with respect to the basis $B$, given below. 
	\begin{equation*}
	\resizebox{1.0 \textwidth}{!} 
	{$
	M(\sigma_{I})=\left( 
		\begin{array}{ccccccccccccccccccc}
		-1& 0& 0& 0& 0& 0& -1& -1& 0& 0& 0& 0& -1& -1& -1& 0& -1& -2& -1 \\ 0& -1& 0& 0& 0& 0& -1& -1& 0& 0& 0& 0& -1& -1& -1& 0& -1& -2& -1 \\ 0& 0& -1& 0& 0& 0& -1& -1& 0& 0& 0& 0& -1& -1& -1& 0& -1& -2& -1 \\ 0& 0& 0& 1& 0& 0& 0& 0& 0& 0& 0& 0& -1& -1& -1& 0& -1& -3& -1 \\ 0& 0& 0& 0& 1& 0& 0& 0& 0& 0& 0& 0& 0& 0& 0& 0& 0& -1& -1 \\ 0& 0& 0& 0& 0& -1& -1& -1& 0& 0& 0& 0& -1& -1& -1& 0& -1& -2& -1 \\ -1& -1& -1& 0& 0& -1& -1& -2& -1& 0& 0& -1& 0& 0& 0& 0& 0& -2& -1 \\ 1& 1& 1& 0& 0& 1& 2& 3& 1& 0& 0& 1& 2& 2& 2& 0& 2& 6& 3 \\ 0& 0& 0& 0& 0& 0& 0& 0& -1& 0& 0& 0& 0& 0& 0& 0& 0& 0& 0 \\ 0& 0& 0& 0& 0& 0& 0& 0& 0& 1& 0& 0& 0& 0& 0& 0& 0& -1& 0 \\ 0& 0& 0& 0& 0& 0& 0& 0& 0& 0& 1& 0& 1& 1& 1& 0& 1& 1& 0 \\ 0& 0& 0& 0& 0& 0& 0& 0& 0& 0& 0& -1& 0& 0& 0& 0& 0& 0& 0 \\ 0& 0& 0& 0& 0& 0& 0& 0& 0& 0& 0& 0& -1& 0& 0& 0& 0& 0& 0 \\ 0& 0& 0& 0& 0& 0& 0& 0& 0& 0& 0& 0& 0& -1& 0& 0& 0& 0& 0 \\ 0& 0& 0& 0& 0& 0& 0& 0& 0& 0& 0& 0& 0& 0& -1& 0& 0& 0& 0 \\ 0& 0& 0& 0& 0& 0& 0& 0& 0& 0& 0& 0& 0& 0& 0& 1& 0& -2& 0 \\ 0& 0& 0& 0& 0& 0& 0& 0& 0& 0& 0& 0& 0& 0& 0& 0& -1& 0& 0 \\ 0& 0& 0& 0& 0& 0& 0& 0& 0& 0& 0& 0& 0& 0& 0& 0& 0& -1& 0 \\ 0& 0& 0& 0& 0& 0& 0& 0& 0& 0& 0& 0& 0& 0& 0& 0& 0& 0& -1
		\end{array}
		\right) 
	$}
	\end{equation*}
	\begin{equation*}
		\resizebox{1.0 \textwidth}{!} 
		{$
		M(\sigma_{2})=\left( 
			\begin{array}{ccccccccccccccccccc}
			0& -1& -1& -1& -1& -1& 0& -2& -1& -1& -1& -1& 0& 0& -1& -1& -1& 1& -1 \\ -1& 0& -1& -1& -1& -1& 0& -2& -1& -1& -1& -1& 0& 0& -1& -1& -1& 1& -1 \\ -1& -1& -2& -1& -1& -1& -1& -3& -1& -1& -1& -1& 0& 0& -1& -1& -1& 1& -1 \\ -1& -1& -1& -2& -1& -1& -1& -3& -1& -1& -1& -1& 0& 0& -1& -1& -1& 2& -1 \\ -1& -1& -1& -1& -2& -1& -1& -3& -1& -1& -1& -1& -1& -1& -1& -1& -1& 0& 0 \\ -1& -1& -1& -1& -1& -2& -1& -3& -1& -1& -1& -1& 0& 0& -1& -1& -1& 1& -2 \\ 0& 0& -1& -1& -1& -1& -1& -2& -1& -1& -1& -1& -1& -1& -1& -1& -1& 1& 0 \\ 2& 2& 3& 3& 3& 3& 2& 7& 3& 3& 3& 3& 1& 1& 3& 3& 3& -3& 2 \\ 0& 0& 0& 0& 0& 0& 0& 0& -1& 0& 0& 0& 0& 0& 0& 0& 0& 0& 0 \\ 0& 0& 0& 0& 0& 0& 0& 0& 0& -1& 0& 0& 0& 0& 0& 0& 0& 1& 0 \\ 0& 0& 0& 0& 0& 0& 0& 0& 0& 0& -1& 0& -1& -1& 0& 0& 0& -1& 1 \\ 0& 0& 0& 0& 0& 0& 0& 0& 0& 0& 0& -1& 0& 0& 0& 0& 0& 0& 1 \\ 0& 0& 0& 0& 0& 0& 0& 0& 0& 0& 0& 0& 1& 0& 0& 0& 0& 0& 0 \\ 0& 0& 0& 0& 0& 0& 0& 0& 0& 0& 0& 0& 0& 1& 0& 0& 0& 0& 0 \\ 0& 0& 0& 0& 0& 0& 0& 0& 0& 0& 0& 0& 0& 0& -1& 0& 0& 0& 0 \\ 0& 0& 0& 0& 0& 0& 0& 0& 0& 0& 0& 0& 0& 0& 0& -1& 0& 2& 0 \\ 0& 0& 0& 0& 0& 0& 0& 0& 0& 0& 0& 0& 0& 0& 0& 0& -1& 0& -2 \\ 0& 0& 0& 0& 0& 0& 0& 0& 0& 0& 0& 0& 0& 0& 0& 0& 0& 1& 0 \\ 0& 0& 0& 0& 0& 0& 0& 0& 0& 0& 0& 0& 0& 0& 0& 0& 0& 0& 1
			\end{array}
			\right) 
		$}
		\end{equation*}
		\begin{equation*}
			\resizebox{1.0 \textwidth}{!} 
			{$
			M(\sigma_{+})=\left( 
				\begin{array}{ccccccccccccccccccc}
				-1& 0& -1& -1& -1& -1& 0& -2& -1& -1& -1& -1& 0& 0& -1& -1& -1& 0& -2 \\ 0& -1& -1& -1& -1& -1& 0& -2& -1& -1& -1& -1& 0& 0& -1& -1& -1& 0& -2 \\ -1& -1& -2& -1& -1& -1& -1& -3& -1& -1& -1& -1& 0& 0& -1& -1& -1& -1& -3 \\ -1& -1& -1& -2& -1& -1& -1& -3& -1& -1& -1& -1& 0& 0& -1& -1& -1& -1& -4 \\ -1& -1& -1& -1& -2& -1& -1& -3& -1& -1& -1& -1& -1& -1& -1& -1& -1& -1& -1 \\ -1& -1& -1& -1& -1& -2& -1& -3& -1& -1& -1& -1& 0& 0& -1& -1& -1& -1& -4 \\ 0& 0& -1& -1& -1& -1& -1& -2& -1& -1& -1& -1& -1& -1& -1& -1& -1& 0& -1 \\ 2& 2& 3& 3& 3& 3& 2& 7& 3& 3& 3& 3& 1& 1& 3& 3& 3& 2& 7 \\ 0& 0& 0& 0& 0& 0& 0& 0& -1& 0& 0& 0& 0& 0& 0& 0& 0& 0& 0 \\ 0& 0& 0& 0& 0& 0& 0& 0& 0& -1& 0& 0& 0& 0& 0& 0& 0& 0& -1 \\ 0& 0& 0& 0& 0& 0& 0& 0& 0& 0& -1& 0& -1& -1& 0& 0& 0& 0& 2 \\ 0& 0& 0& 0& 0& 0& 0& 0& 0& 0& 0& -1& 0& 0& 0& 0& 0& 0& 1 \\ 0& 0& 0& 0& 0& 0& 0& 0& 0& 0& 0& 0& 0& 1& 0& 0& 0& 0& 0 \\ 0& 0& 0& 0& 0& 0& 0& 0& 0& 0& 0& 0& 1& 0& 0& 0& 0& 0& 0 \\ 0& 0& 0& 0& 0& 0& 0& 0& 0& 0& 0& 0& 0& 0& -1& 0& 0& 0& 0 \\ 0& 0& 0& 0& 0& 0& 0& 0& 0& 0& 0& 0& 0& 0& 0& -1& 0& 0& -2 \\ 0& 0& 0& 0& 0& 0& 0& 0& 0& 0& 0& 0& 0& 0& 0& 0& -1& 0& -2 \\ 0& 0& 0& 0& 0& 0& 0& 0& 0& 0& 0& 0& 0& 0& 0& 0& 0& -1& -2 \\ 0& 0& 0& 0& 0& 0& 0& 0& 0& 0& 0& 0& 0& 0& 0& 0& 0& 0& 1
				\end{array}
				\right) 
			$}
			\end{equation*}
			\begin{equation*}
				\resizebox{1.0 \textwidth}{!} 
				{$
				M(\sigma_{-})=\left( 
					\begin{array}{ccccccccccccccccccc}
					-1& -2& -1& -1& -1& -1& -1& -3& -1& -1& -1& -1& -1& -1& -1& -1& -1& 0& 0 \\ -2& -1& -1& -1& -1& -1& -1& -3& -1& -1& -1& -1& -1& -1& -1& -1& -1& 0& 0 \\ -1& -1& -2& -1& -1& -1& -1& -3& -1& -1& -1& -1& -1& -1& -1& -1& -1& 1& 1 \\ -1& -1& -1& -2& -1& -1& -1& -3& -1& -1& -1& -1& -1& -1& -1& -1& -1& 2& 2 \\ -1& -1& -1& -1& -2& -1& -1& -3& -1& -1& -1& -1& -1& -1& -1& -1& -1& 0& 0 \\ -1& -1& -1& -1& -1& -2& -1& -3& -1& -1& -1& -1& -1& -1& -1& -1& -1& 1& 1 \\ -1& -1& -1& -1& -1& -1& -2& -3& -1& -1& -1& -1& -1& -1& -1& -1& -1& 0& 0 \\ 3& 3& 3& 3& 3& 3& 3& 8& 3& 3& 3& 3& 3& 3& 3& 3& 3& -2& -2 \\ 0& 0& 0& 0& 0& 0& 0& 0& -1& 0& 0& 0& 0& 0& 0& 0& 0& 0& 0 \\ 0& 0& 0& 0& 0& 0& 0& 0& 0& -1& 0& 0& 0& 0& 0& 0& 0& 1& 1 \\ 0& 0& 0& 0& 0& 0& 0& 0& 0& 0& -1& 0& 0& 0& 0& 0& 0& -1& -1 \\ 0& 0& 0& 0& 0& 0& 0& 0& 0& 0& 0& -1& 0& 0& 0& 0& 0& 0& 0 \\ 0& 0& 0& 0& 0& 0& 0& 0& 0& 0& 0& 0& 0& -1& 0& 0& 0& 0& 0 \\ 0& 0& 0& 0& 0& 0& 0& 0& 0& 0& 0& 0& -1& 0& 0& 0& 0& 0& 0 \\ 0& 0& 0& 0& 0& 0& 0& 0& 0& 0& 0& 0& 0& 0& -1& 0& 0& 0& 0 \\ 0& 0& 0& 0& 0& 0& 0& 0& 0& 0& 0& 0& 0& 0& 0& -1& 0& 2& 2 \\ 0& 0& 0& 0& 0& 0& 0& 0& 0& 0& 0& 0& 0& 0& 0& 0& -1& 0& 0 \\ 0& 0& 0& 0& 0& 0& 0& 0& 0& 0& 0& 0& 0& 0& 0& 0& 0& 1& 2 \\ 0& 0& 0& 0& 0& 0& 0& 0& 0& 0& 0& 0& 0& 0& 0& 0& 0& 0& -1
					\end{array}
					\right) 
				$}
				\end{equation*}

	\end{Rmk}

\bibliographystyle{amsalpha}

	\newcommand{\etalchar}[1]{$^{#1}$}
	\providecommand{\bysame}{\leavevmode\hbox to3em{\hrulefill}\thinspace}
	\providecommand{\MR}{\relax\ifhmode\unskip\space\fi MR }
	\providecommand{\MRhref}[2]{%
		\href{http://www.ams.org/mathscinet-getitem?mr=#1}{#2}
	}
	\providecommand{\href}[2]{#2}

\end{document}